\documentclass[12pt]{article}
\usepackage{amssymb}
\usepackage{amsfonts}
\usepackage{amsmath}
\usepackage{amsbsy}
\usepackage[small,nohug,heads=littlevee]{diagrams}
\diagramstyle[labelstyle=\scriptstyle]
\newcommand{\be}{\begin{equation}}
\newcommand{\ee}{\end{equation}}
\newcommand{\bea}{\begin{eqnarray}}
\newcommand{\eea}{\end{eqnarray}}
\newcommand{\ba}{\begin{array}}
\newcommand{\ea}{\end{array}}

\newcommand{\bc}{\begin{center}}
\newcommand{\ec}{\end{center}}
\newcommand{\ben}{\begin{enumerate}}
\newcommand{\een}{\end{enumerate}}
\newcommand{\bfi}{\begin{figure}}
\newcommand{\efi}{\end{figure}}

\newcommand{\bq}{\begin{quote}}
\newcommand{\eq}{\end{quote}}
\newcommand{\bqu}{\begin{quotation}}
\newcommand{\equ}{\end{quotation}}
\newenvironment{emphit}{\begin{itemize}}{\end{itemize}}
\newcommand{\bemp}{\begin{emphit}}
\newcommand{\eemp}{\end{emphit}}

\newcommand{\bt}{\begin{tabular}}
\newcommand{\et}{\end{tabular}}

\newtheorem{myth}{Theorem}[section]
\newtheorem{mylem}{Lemma}[section]

\newtheorem{mydef}{Definition}[section]

\newtheorem{myexam}{Example}[section]

\begin{document}
\date{}
\title{On Another Two Cryptographic Identities In Universal Osborn Loops 
}
\author{T. G. Jaiy\'e\d ol\'a\thanks{All correspondence to be addressed to this author.} \\
Department of Mathematics,\\
Obafemi Awolowo University,\\
Ile Ife 220005, Nigeria.\\
jaiyeolatemitope@yahoo.com\\tjayeola@oauife.edu.ng \and
J. O. Ad\'en\'iran \\
Department of Mathematics,\\
University of Agriculture, \\
Abeokuta 110101, Nigeria.\\
ekenedilichineke@yahoo.com\\
adeniranoj@unaab.edu.ng} \maketitle
\begin{abstract}
In this study, by establishing an identity for universal Osborn
loops, two other identities(of degrees $4$ and $6$) are deduced from it and they are
recognized and recommended for cryptography in a similar spirit in
which the cross inverse property(of degree $2$) has been used by Keedwell following
the fact that it was observed that universal Osborn loops that do
not have the 3-power associative property or weaker forms of;
inverse property, power associativity and diassociativity to mention
a few, will have cycles(even long ones). These identities are found
to be cryptographic in nature for universal Osborn loops and thereby
called cryptographic identities. They were also found applicable to security patterns, arrangements and networks which the CIP may not be applicable to.
\end{abstract}

\section{Introduction}
\paragraph{}
Let $L$ be a non-empty set. Define a binary operation ($\cdot $) on
$L$ : If $x\cdot y\in L$ for all $x, y\in L$, $(L, \cdot )$ is
called a groupoid. If the system of equations ;
\begin{displaymath}
a\cdot x=b\qquad\textrm{and}\qquad y\cdot a=b
\end{displaymath}
have unique solutions for $x$ and $y$ respectively, then $(L, \cdot
)$ is called a quasigroup. Furthermore, if there exists a unique
element $e\in L$ called the identity element such that for all $x\in
L$, $x\cdot e=e\cdot x=x$, $(L, \cdot )$ is called a loop. We write
$xy$ instead of $x\cdot y$, and stipulate that $\cdot$ has lower
priority than juxtaposition among factors to be multiplied. For
instance, $x\cdot yz$ stands for x(yz). For each $x\in L$, the
elements $x^\rho =xJ_\rho ,x^\lambda =xJ_\lambda\in L$ such that
$xx^\rho=e=x^\lambda x$ are called the right, left inverses of $x$
respectively. $x^{\lambda^i}=(x^\lambda )^\lambda$ and
$x^{\rho^i}=(x^\rho )^\rho$ for $i\ge 1$.
\begin{mydef}\label{0:1}
A loop $(G,\cdot ,/,\backslash ,e)$ is a set $G$ together with three
binary operations ($\cdot $), ($/$), ($\backslash$) and one nullary
operation $e$ such that
\begin{description}
\item[(i)] $x\cdot (x\backslash y)=y$, $(y/x)\cdot x=y$ for all
$x,y\in G$,
\item[(ii)] $x\backslash (x\cdot y)=y$, $(y\cdot x)/x=y$ for all
$x,y\in G$ and
\item[(iii)] $x\backslash x=y/y$ or $e\cdot x=x$ for all
$x,y\in G$.
\end{description}
\end{mydef}
We also stipulate that ($/$) and ($\backslash$) have higher priority
than ($\cdot $) among factors to be multiplied. For instance,
$x\cdot y/z$ and $x\cdot y\backslash z$ stand for $x(y/z)$ and
$x\cdot (y\backslash z)$ respectively.

The left and right translation maps of $G$, $L_x$ and $R_x$
respectively can be defined by
\begin{displaymath}
yL_x=x\cdot y\qquad\textrm{and}\qquad yR_x=y\cdot x.
\end{displaymath}
Let
\begin{displaymath}
x\backslash y =yL_x^{-1}=y\mathbb{L}_x\qquad\textrm{and}\qquad
x/y=xR_y^{-1}=x\mathbb{R}_y.
\end{displaymath}

$L$ is called a weak inverse property loop (WIPL) if and only if it
obeys the weak inverse property (WIP);
\begin{displaymath}
xy\cdot z=e~\textrm{implies}~x\cdot yz=e~\textrm{for all}~x,y,z\in L
\end{displaymath}
while $L$ is called a cross inverse property loop (CIPL) if and only
if it obeys the cross inverse property (CIP);
\begin{displaymath}
xy\cdot x^\rho=y.
\end{displaymath}
The triple $\alpha =(A,B,C)$ of bijections on a loop $(L,\cdot )$ is
called an autotopism of the loop if and only if
\begin{displaymath}
xA\cdot yB=(x\cdot y)C~\textrm{for all}~x,y\in L.
\end{displaymath}
Such triples form a group $AUT(L,\cdot )$ called the autotopism
group of $(L,\cdot )$. In case the three bijections are the same i.e
$A=B=C$, then any of them is called an automorphism and the group
$AUM(L,\cdot )$ which such forms is called the automorphism group of
$(L,\cdot )$. For an overview of the theory of loops, readers may
check \cite{phd3,phd41,phd39,phd49,phd42,phd75}.

Osborn \cite{phd89}, while investigating the universality of WIPLs
discovered that a universal WIPL $(G,\cdot )$ obeys the identity
\begin{equation}\label{eq:1}
yx\cdot (z\theta_y\cdot y)=(y\cdot xz)\cdot y~\textrm{for all}~x,y,z\in G
\end{equation}
\begin{displaymath}
\textrm{where}~\theta_y=L_yL_{y^\lambda}=R_{y^\rho}^{-1}R_y^{-1}=L_yR_yL_y^{-1}R_y^{-1}.
\end{displaymath}
A loop that necessarily and sufficiently satisfies this identity is
called an Osborn loop.

\paragraph{}
Eight years after Osborn's \cite{phd89} 1960 work on WIPL, in 1968,
Huthnance  Jr. \cite{phd44} studied the theory of generalized
Moufang loops. He named a loop that obeys (\ref{eq:1}) a generalized
Moufang loop and later on in the same thesis, he called them
M-loops. On the other hand, he called a universal WIPL an Osborn
loop and this same definition was adopted by Chiboka \cite{phd96}.
Basarab \cite{phd148,phd46,phd137} and Basarab and Belioglo
\cite{phd170} dubbed a loop $(G,\cdot )$ satisfying any of the
following equivalent identities an Osborn loop:
\begin{equation}
OS_2~:~x(yz\cdot x)=(x^\lambda\backslash y)\cdot  zx
\end{equation}
\begin{equation}
OS_3~:~(x\cdot yz)x=xy\cdot (zE_x^{-1}\cdot x)
\end{equation}
\begin{displaymath}
\textrm{where}~E_x=R_xR_{x^\rho}=(L_xL_{x^\lambda})^{-1}=R_xL_xR_x^{-1}L_x^{-1}~\textrm{for
all}~x,y,z\in G
\end{displaymath} and the binary operations '$\backslash$' and
'$/$' are respectively defines as ; $z=x\cdot y$ if and only if
$x\backslash z=y$ if and only if $z/y=x$ for all $x,y,z\in G$.

It will look confusing if both Basarab's and Huthnance's definitions
of an Osborn loop are both adopted because an Osborn loop of Basarab
is not necessarily a universal WIPL(Osborn loop of Huthnance). So in
this work, Huthnance's definition of an Osborn loop will be dropped
while we shall stick to that of Basarab which was actually adopted
by M. K. Kinyon \cite{phd33} who revived the study of Osborn loops
in 2005 at a conference tagged "Milehigh Conference on Loops,
Quasigroups and Non-associative Systems" held at the University of
Denver, where he presented a talk titled "A Survey of Osborn Loops".

Let $t=x^\lambda\backslash y$ in $OS_2$, then $y=x^\lambda t$ so
that we now have an equivalent identity
$$x[(x^\lambda y)z\cdot x]=y\cdot  zx.$$
Huthnance \cite{phd44} was able to deduce some properties of $E_x$
relative to (\ref{eq:1}). $E_x=E_{x^\lambda}=E_{x^\rho}$. So, since
$E_x=R_xR_{x^\rho}$, then $E_x=E_{x^\lambda}=R_{x^\lambda}R_{x}$ and
$E_x=(L_{x^\rho}L_x)^{-1}$. So, we now have the following equivalent
identity defining an Osborn loop.
\begin{equation}\label{eq:1.1}
\textrm{OS$_0$}~:~x(yz\cdot x)=x(yx^\lambda \cdot x)\cdot zx
\end{equation}
\begin{mydef}
A loop $(Q,\cdot )$ is called:
\begin{description}
\item[(a)] a 3 power associative property loop(3-PAPL) if and only if $xx\cdot
x=x\cdot xx$ for all $x\in Q$.
\item[(b)] a left self inverse property loop(LSIPL) if and only if
$x^\lambda\cdot xx=x$ for all $x\in Q$.
\item[(c)] a right self inverse property loop(RSIPL) if and only if
$xx\cdot x^\rho =x$ for all $x\in Q$.
\end{description}
\end{mydef}

The identities describing the most popularly known varieties of
Osborn loops are given below.
\begin{mydef}
A loop $(Q,\cdot )$ is called:
\begin{description}
\item[(a)] a VD-loop if and only if
\begin{displaymath}
(\cdot )_x=(\cdot )^{L_x^{-1}R_x}\qquad\textrm{and}\qquad {}_x(\cdot
)=(\cdot )^{R_x^{-1}L_x} \end{displaymath} i.e $R_x^{-1}L_x\in
PS_\lambda (Q,\cdot )$ with companion $c=x$ and $L_x^{-1}R_x\in
PS_\rho (Q,\cdot )$ with companion $c=x$ for all $x\in Q$ where
$PS_\lambda (Q,\cdot )$ and $PS_\rho (Q,\cdot )$ are respectively
the left and right pseudo-automorphism groups of $Q$.\qquad Basarab
\cite{phd137}
\item[(b)] a Moufang loop if and only if the identity
\begin{displaymath}
(xy)\cdot (zx)=(x\cdot yz)x
\end{displaymath}
holds in $Q$.
\item[(c)] a conjugacy closed loop(CC-loop) if and only if the
identities
\begin{displaymath}
x\cdot yz=(xy)/x\cdot xz\qquad\textrm{and}\qquad zy\cdot x=zx\cdot
x\backslash(yx)
\end{displaymath}
hold in $Q$.
\item[(d)] a universal WIPL if and only if the identity
\begin{displaymath}
x(yx)^\rho=y^\rho\qquad\textrm{or}\qquad(xy)^\lambda x=y^\lambda
\end{displaymath}
holds in $Q$ and all its isotopes.
\end{description}
\end{mydef}
All these three varieties of Osborn loops and universal WIPLs are
universal Osborn loops. CC-loops and VD-loops are G-loops. G-loops
are loops that are isomorphic to all their loop isotopes. Kunen
\cite{phd185} has studied them.

In the multiplication group  $\mathcal{M}\textrm{ult}(Q)$ of a loop
$(G,\cdot )$ are found three important permutations, namely, the
right, left and middle inner mappings $R_{(x,y)}=R_xR_yR_{xy}^{-1}$,
$L_{(x,y)}=L_xL_yL_{yx}^{-1}$ and $T_{(x)}=R_xL_x^{-1}$ respectively
which form the right inner mapping group $\textrm{Inn}_\lambda(G)$,
left inner mapping group $\textrm{Inn}_\rho (G)$ and the middle
inner mapping $\textrm{Inn}_\mu (G)$. In a Moufang loop $G$,
$R_{(x,y)},L_{(x,y)},T_{(x)}\in PS_\rho(G)$ with companions
$(x,y),(x^{-1},y^{-1}),x^{-3}\in G$ respectively.
\begin{myth}\label{1:3.02}(Kinyon \cite{phd33})

Let $G$ be an Osborn loop. $R_{(x,y)}\in PS_\rho(G)$ with companion
$(xy)^\lambda (y^\lambda\backslash x)$ and $L_{(x,y)}\in
PS_\lambda(G)~\forall~x,y\in G$. Furthermore,
$R_{(x,y)}^{-1}=[L_{y^\rho}^{-1},R_x^{-1}]=L_{(y^\lambda
,x^\lambda)}~\forall~x,y\in G$.
\end{myth}
The second part of Theorem~\ref{1:3.02} is trivial for Moufang
loops. For CC-loops, it was first observed by Dr\'apal and then
later by Kinyon and Kunen \cite{phd47}.
\begin{myth}\label{1:3.05}
Let $G$ be an Osborn loop.
$\textrm{Inn}_\rho(G)=\textrm{Inn}_\lambda(G)$.
\end{myth}
Still mysterious are the middle inner mappings $T_{(x)}$ of an
Osborn loop. In a  Moufang loop, $T_{(x)}\in PS_\rho$ with a
companion $x^{-3}$ while in a CC-loop, $T_{(x)}\in PS_\lambda$ with
companion $x$. So, Kinyon \cite{phd33} possessed a question asking of which group(whether $PS_\rho$ and $PS_\lambda$) to which $T_{(x)}$ belongs to in case of an arbitrary Osborn loop and what its companion will be.
\begin{myth}\label{1:3.050}(Kinyon \cite{phd33})

In an Osborn loop $G$ with centrum $C(G)$ and center $Z(G)$:
\begin{enumerate}
\item If $T_{(a)}\in AUM(G)$, then $a\cdot aa=aa\cdot a\in N(G)$. Thus, for all $a\in C(G)$, $a^3\in Z(G)$.
\item If $(xx)^\rho =x^\rho x^\rho$ holds, then
$x^{\rho\rho\rho\rho\rho\rho}=x$ for all $x\in G$.
\end{enumerate}
\end{myth}

Some basic loop properties such as flexibility, left alternative
property(LAP), left inverse property(LIP), right alternative
property(RAP), right inverse property(RIP), anti-automorphic inverse
property(AAIP) and the cross inverse property(CIP) have been found
to force an Osborn loop to be a Moufang loop. This makes the study
of Osborn loops more challenging and care must be taking not to
assume any of these properties at any point in time except the WIP,
automorphic inverse property and some other generalizations of the
earlier mentioned loop properties(LAP, LIP, e.t.c.).

\begin{mylem}\label{1:4.1}
An Osborn loop that is flexible or which has the LAP or RAP or LIP
or RIP or AAIP is a Moufang loop. But an Osborn loop that is
commutative or which has the CIP is a commutative Moufang loop.
\end{mylem}
\begin{myth}\label{0:2}(Basarab, \cite{phd46})

If an Osborn loop is of exponent 2, then it is an abelian group.
\end{myth}

\begin{myth}\label{1:3.13}(Huthnance \cite{phd44})

Let $G$ be a WIPL. $G$ is a universal WIPL if and only if $G$ is an
Osborn loop.
\end{myth}

\begin{mylem}\label{3:13}(Lemma~2.10, Huthnance \cite{phd44})

Let $L$ be a WIP Osborn loop. If $a=x^\rho x$, then for all $x\in
L$:
\begin{displaymath}
xa=x^{\lambda^2},~ax^\lambda =x^\rho,~x^\rho
a=x^\lambda,~ax=x^{\rho^2},~xa^{-1}=ax,~a^{-1}x^\lambda =x^\lambda
a,~a^{-1}x^\rho =x^\rho a.
\end{displaymath}
or equivalently
\begin{displaymath}
J_\lambda~:~x\mapsto~x\cdot x^\rho x,~J_\rho~:~x\mapsto~x^\rho
x\cdot x^\lambda,~J_\lambda~:~x\mapsto~x^\rho\cdot x^\rho x
,~J_\rho^2~:~x\mapsto~x^\rho x\cdot x,
\end{displaymath}
\begin{displaymath}
x(x^\rho x)^{-1}=(x^\rho x)x,~(x^\rho x)^{-1}x^\lambda =x^\lambda
\cdot x^\rho x,~(x^\rho x)^{-1}x^\rho =x^\rho (x^\rho x).
\end{displaymath}
\end{mylem}

Consider $(G,\cdot )$ and $(H,\circ )$ been two distinct
groupoids or quasigroups or loops. Let $A,B$ and $C$ be three bijective
mappings, that map $G$ onto $H$. The triple $\alpha =(A,B,C)$ is
called an isotopism of $(G,\cdot )$ onto $(H,\circ )$ if and only if
$$xA\circ yB=(x\cdot y)C~\forall~x,y\in G.$$
So, $(H,\circ )$ is called a groupoid(quasigroup, loop) isotope of
$(G,\cdot )$.

If $C=I$ is the identity map on $G$ so that $H=G$, then the triple
$\alpha =(A,B,I)$ is called a principal isotopism of $(G,\cdot )$
onto $(G,\circ )$ and $(G,\circ)$ is called a principal isotope of
$(G,\cdot )$. Eventually, the equation of relationship now becomes
$$x\cdot y=xA\circ yB~\forall~x,y\in G$$
which is easier to work with. But if $A=R_g$ and $B=L_f$, for some
$f,g\in G$, the relationship now becomes
$$x\cdot y=xR_g\circ yL_f~\forall~x,y\in G$$
or
$$x\circ y=xR_g^{-1}\cdot yL_f^{-1}~\forall~x,y\in G.$$
With this new form, the triple $\alpha =(R_g,L_f,I)$ is called an
$f,g$-principal isotopism of $(G,\cdot )$ onto $(G,\circ )$, $f$ and
$g$ are called translation elements of $G$ or at times written in
the pair form $(g,f)$, while $(G,\circ )$ is called an
$f,g$-principal isotope of $(G,\cdot )$.

The last form of $\alpha$ above gave rise to an important result in
the study of loop isotopes of loops.

\begin{myth}\label{1:1}(Bruck \cite{phd41})

Let $(G,\cdot )$ and $(H,\circ )$ be two distinct isotopic loops.
For some $f,g\in G$, there exists an $f,g$-principal isotope
$(G,\ast )$ of $ (G,\cdot )$ such that $(H,\circ )\cong (G,\ast )$.
\end{myth}

With this result, to investigate the isotopic invariance of an
isomorphic invariant property in loops, one simply needs only to
check if the property in consideration is true in all
$f,g$-principal isotopes of the loop. A property is isotopic
invariant if whenever it holds in the domain loop i.e $(G,\cdot )$
then it must hold in the co-domain loop i.e $(H,\circ )$ which is an
isotope of the formal. In such a situation, the property in
consideration is said to be a universal property hence the loop is
called a universal loop relative to the property in consideration as
often used by Nagy and Strambach \cite{phd88} in their algebraic and
geometric study of the universality of some types of loops. For
instance, if every loop isotope of a loop with property $\mathcal{P}$ also has the property $\mathcal{P}$, then the formal is called a universal $\mathcal{P}$ loop. So, we can now restate Theorem~\ref{1:1} as :

\begin{myth}\label{1:2}
Let $(G,\cdot )$ be a loop with an
isomorphic invariant property $\mathcal{P}$. $(G,\cdot )$ is a universal $\mathcal{P}$ loop if and only if every $f,g$-principal isotope $(G,\ast )$ of $(G,\cdot )$ has the $\mathcal{P}$ property.
\end{myth}

\begin{mydef}(Universal Osborn Loop)
A loop is called a universal Osborn loop if all its loop isotopes
are Osborn loops.
\end{mydef}

The aim of this study is to identify some identities that are
appropriate for cryptography in universal Osborn loops. These
identities hold in universal Osborn loops like CC-loops, introduced
by Goodaire and Robinson \cite{phd91,phd48}, whose algebraic
structures have been studied by Kunen \cite{phd78} and some recent
works of Kinyon and Kunen \cite{phd36,phd47}, Phillips et. al.
\cite{phd35}, Dr\'apal \cite{phd37,phd38,phd98,phd107}, Cs\"org\H o
et. al. \cite{phd106,phd104,phd108} and VD-loops whose study is yet
to be explored. In this study, by establishing an identity for universal Osborn
loops, two other identities(of degrees $4$ and $6$) are deduced from it and they are
recognized and recommended for cryptography in a similar spirit in
which the cross inverse property(of degree $2$) has been used by Keedwell following
the fact that it was observed that universal Osborn loops that do
not have the 3-power associative property or weaker forms of;
inverse property, power associativity and diassociativity to mention
a few, will have cycles(even long ones). These identities are found
to be cryptographic in nature for universal Osborn loops and thereby
called cryptographic identities. They were also found applicable to security patterns, arrangements and networks which the CIP may not be applicable to.

We shall make use of the following results.
\paragraph{Results of Bryant and Schneider \cite{phd92}}
\begin{myth}\label{0:3}
Let $(Q, \cdot ,\backslash ,/)$ be a quasigroup. If $Q(a,b,\circ
)\overset{\theta}{\cong} Q(c,d,\ast )$ for any $a,b,c,d\in Q$, then $Q(f,g,\circledcirc
)\overset{\theta}{\cong} Q\big((f\cdot b)\theta/d, c\backslash
(a\cdot g)\theta ,\star\big)$ for any $a,b,c,d,f,g\in Q$. If $(Q, \cdot )$ is a loop, then
\begin{displaymath}
(f\cdot b)\theta/d=[f\cdot (a\backslash
c\theta^{-1})]\theta~\textrm{and}~c\backslash (a\cdot
g)\theta =[(d\theta^{-1}/b)\cdot g]\theta~\textrm{for any}~a,b,c,d,f,g\in Q.
\end{displaymath}
\end{myth}

\section{Main Results}
\subsection{Identities In Universal Osborn Loops}
\begin{myth}\label{1:4}
Let $(Q, \cdot ,\backslash ,/)$ be an Osborn loop, $(Q,\circ )$ an arbitrary
principal isotope of $(Q,\cdot )$ and $(Q,\ast )$ some principal
isotopes of $(Q,\cdot )$. Let $\phi (x,u,v)=(u\backslash
([(uv)/(u\backslash (xv))]v))$ and $\gamma=\gamma
(x,u,v)=\mathbb{R}_vR_{[u\backslash (xv)]}\mathbb{L}_uL_x$ for all
$x,u,v\in Q$, then $(Q, \cdot ,\backslash ,/)$ is a universal Osborn
loop if and only if the commutative diagram
\begin{equation}\label{eq:7}
\begin{diagram}
&                                                                               &(Q,\ast )\\
&\ruTo^{(R_{\phi (x,u,v)},L_u,I)}_{}&\dTo^{(\gamma,\gamma,\gamma)}_{\textrm{isomorphism}}\\
(Q,\cdot ) &\rTo^{(R_v,L_x,I)}_{\textrm{principal isotopism}}
&(Q,\circ)
\end{diagram}
\end{equation}
holds.
\end{myth}
{\bf Proof}\\
Let $\mathcal{Q}=(Q, \cdot ,\backslash ,/)$ be an Osborn loop with
any arbitrary principal isotope $\mathfrak{Q}=(Q, \blacktriangle
,\nwarrow ,\nearrow )$ such that
\begin{displaymath}
x\blacktriangle y=xR_v^{-1}\cdot yL_u^{-1}=(x/v)\cdot (u\backslash
y)~\forall~u,v\in Q.
\end{displaymath}
If $\mathcal{Q}$ is a universal Osborn loop then, $\mathfrak{Q}$ is
an Osborn loop. $\mathfrak{Q}$ obeys identity OS$_0$ implies
\begin{equation}\label{eq:3}
x\blacktriangle [(y\blacktriangle z)\blacktriangle
x]=\{x\blacktriangle [(y\blacktriangle x^{\lambda '})\blacktriangle
x]\}\blacktriangle (z\blacktriangle x)
\end{equation}
where $x^{\lambda '}=xJ_{\lambda '}$ is the left inverse of $x$ in
$\mathfrak{Q}$. The identity element of the loop $\mathfrak{Q}$ is
$uv$. So,
\begin{displaymath}
x\blacktriangle y=xR_v^{-1}\cdot yL_u^{-1}~\textrm{implies}~
y^{\lambda '}\blacktriangle y=y^{\lambda '}R_v^{-1}\cdot
yL_u^{-1}=uv~\textrm{implies}
\end{displaymath}
\begin{displaymath}
y^{\lambda '}R_v^{-1}R_{yL_u^{-1}}=uv~\textrm{implies}~yJ_{\lambda
'}=(uv)R_{yL_u^{-1}}^{-1}R_v=(uv)R_{(u\backslash
y)}^{-1}R_v=[(uv)/(u\backslash y)]v.
\end{displaymath}
Thus, using the fact that
\begin{displaymath}
x\blacktriangle y=(x/v)\cdot (u\backslash y),
\end{displaymath}
$\mathfrak{Q}$ is an Osborn loop if and only if{\footnotesize
\begin{displaymath}
(x/v)\cdot u\backslash\{[(y/v)\cdot (u\backslash z)]/v\cdot
(u\backslash x)\}=((x/v)\cdot
u\backslash\{[(y/v)(u\backslash([(uv)/(u\backslash x)]v))]/v\cdot
(u\backslash x)\})/v\cdot u\backslash[(z/v)(u\backslash x)].
\end{displaymath}}
Do the following replacements:
\begin{displaymath}
x'=x/v\Rightarrow x=x'v,~z'=u\backslash z\Rightarrow
z=uz',~y'=y/v\Rightarrow y=y'v
\end{displaymath}
we have{\footnotesize
\begin{displaymath}
x'\cdot u\backslash\{(y'z')/v\cdot [u\backslash (x'v)]\}=(x'\cdot
u\backslash\{[y'(u\backslash([(uv)/(u\backslash (x'v))]v))]/v\cdot
[u\backslash (x'v)]\})/v  \cdot u\backslash[((uz')/v)(u\backslash
(x'v))].
\end{displaymath}}
This is precisely identity OS$_0'$ below by replacing $x'$, $y'$ and
$z'$ by $x$, $y$ and $z$ respectively. {\footnotesize
\begin{displaymath}
\underbrace{ x\cdot u\backslash\{(yz)/v\cdot [u\backslash
(xv)]\}=(x\cdot u\backslash\{[y(u\backslash([(uv)/(u\backslash
(xv))]v))]/v\cdot [u\backslash (xv)]\})/v  \cdot
u\backslash[((uz)/v)(u\backslash (xv))].}_{\textrm{OS$_0'$}}
\end{displaymath}}
Writing identity OS$_0'$ in autotopic form, we will obtain the fact
that the triple $\big(\alpha (x,u,v),\beta (x,u,v),\gamma
(x,u,v)\big)\in AUT(\mathcal{Q})$ for all $x,u,v\in Q$ where $\alpha
(x,u,v)=R_{_{(u\backslash([(uv)/(u\backslash
(xv))]v))}}\mathbb{R}_vR_{[u\backslash
(xv)]}\mathbb{L}_uL_x\mathbb{R}_v,~\beta
(x,u,v)=L_u\mathbb{R}_vR_{[u\backslash (xv)]}\mathbb{L}_u$ and
$\gamma (x,u,v)=\mathbb{R}_vR_{[u\backslash (xv)]}\mathbb{L}_uL_x$
are elements of $\mathcal{M}\textrm{ult}(Q)$. The triple
\begin{displaymath}
\big(\alpha (x,u,v),\beta (x,u,v),\gamma
(x,u,v)\big)=\Big(R_{_{(u\backslash([(uv)/(u\backslash
(xv))]v))}}\gamma\mathbb{R}_v,L_u\gamma\mathbb{L}_x,\gamma\Big)
\end{displaymath}
can be written as the following compositions
$\Big(R_{_{(u\backslash([(uv)/(u\backslash
(xv))]v))}},L_u,I\Big)(\gamma,\gamma,\gamma)(\mathbb{R}_v,\mathbb{L}_x,I)$.
Let $(Q,\circ )$ be an arbitrary principal isotope of $(Q,\cdot )$
and $(Q,\ast )$ a particular principal isotope of $(Q,\cdot )$ under the isotopism $(R_{\phi
(x,u,v)},L_u,I)$ where $\phi (x,u,v)=(u\backslash ([(uv)/(u\backslash (xv))]v))$. Then, the
composition above can be expressed as:
\begin{displaymath}
(Q,\cdot )\xrightarrow[\textrm{principal isotopism}]{(R_{\phi
(x,u,v)},L_u,I)}(Q,\ast
)\xrightarrow[\textrm{isomorphism}]{(\gamma,\gamma,\gamma)}(Q,\circ
)\xrightarrow[\textrm{principal
isotopism}]{(\mathbb{R}_v,\mathbb{L}_x,I)}(Q,\cdot ).
\end{displaymath}

The proof of the converse is as follows. Let $\mathcal{Q}=(Q, \cdot
,\backslash ,/)$ be an Osborn loop. Assuming that the composition in
equation (\ref{eq:7}) holds, then doing the reverse of the proof of
necessity, $\big(\alpha (x,u,v),\beta (x,u,v),\gamma (x,u,v)\big)\in
AUT(\mathcal{Q})$ for all $x,u,v\in Q$ which means that
$\mathcal{Q}$ obeys identity OS$_0'$ hence, it will be observed that
equation (\ref{eq:3}) is true for any arbitrary $u,v$-principal
isotope $\mathfrak{Q}=(Q, \blacktriangle ,\nwarrow ,\nearrow )$ of
$\mathcal{Q}$. So, every $f,g$-principal isotope $\mathfrak{Q}$ of
$\mathcal{Q}$ is an Osborn loop. Following Theorem~\ref{1:2},
$\mathcal{Q}$ is a universal Osborn loop if and only if
$\mathfrak{Q}$ is an Osborn loop.

\begin{myth}\label{1:8}
A universal Osborn loop $(Q, \cdot ,\backslash ,/)$ obeys the
identity {\footnotesize
\begin{displaymath}
\underbrace{(u[x\backslash (zv)])/[u\backslash (xv)]\cdot
v=(\{u\cdot x\backslash\{z\cdot u\backslash ((u/v)[u\backslash
(xv)])\}\}/[u\backslash (xv)]\cdot v)\cdot u\backslash
([(uv)/(u\backslash (xv))]v).}_{\textrm{OSI$_0^1$}}
\end{displaymath}}
for all $x,z,u,v\in Q$.
\begin{align*}
\textrm{Furthermore,}~\underbrace{z=x\cdot \{[x\backslash
(zx)]/x\cdot
x^\lambda\}x}_{\textrm{OSI$_0^{1.1}$}}\qquad\textrm{and}\qquad\underbrace{(x^\lambda\cdot
xy)x^\lambda\cdot x=y}_{\textrm{double left inverse property(DLIP)}}
\end{align*}
are also satisfied for all $x,y,z\in Q$.
\end{myth}
{\bf Proof}\\
By equation (\ref{eq:7}) of Theorem~\ref{1:4}, it can be deduced
that if $(Q,\circ )$ and $(Q,\ast )$ are principal isotopes of
$(Q,\cdot )$ and $\gamma (x,u,v)=\mathbb{R}_vR_{[u\backslash
(xv)]}\mathbb{L}_uL_x$, then
\begin{displaymath}
(Q,x,v,\circ )\overset{~~\gamma^{-1}}{\cong}(Q,u,\phi (x,u,v),\ast
)~\textrm{where}~\phi (x,u,v)=(u\backslash ([(uv)/(u\backslash
(xv))]v))~\textrm{for all}~x,u,v\in Q.
\end{displaymath}
Let $Q(z,y,\circledcirc )$ be an arbitrary principal isotope of
$(Q,\cdot )$. We now switch to Theorem~\ref{0:3}. Let $a=x$, $b=v$,
$c=u$, $d=\phi (x,u,v)=(u\backslash ([(uv)/(u\backslash (xv))]v))$,
$f=z$ and $g=y$. $\theta =\gamma
(x,u,v)^{-1}=\mathbb{L}_xL_u\mathbb{R}_{[u\backslash (xv)]}R_v$
while $\theta^{-1}=\gamma (x,u,v)=\mathbb{R}_vR_{[u\backslash
(xv)]}\mathbb{L}_uL_x$.
\begin{displaymath}
(f\cdot b)\theta/d=\{(u[x\backslash (zv)])/[u\backslash (xv)]\cdot
v\}/\{u\backslash ([(uv)/(u\backslash (xv))]v)\}~\textrm{and}
\end{displaymath}
\begin{displaymath}
[f\cdot (a\backslash c\theta^{-1})]\theta =\{u\cdot
x\backslash\{z\cdot u\backslash ((u/v)[u\backslash
(xv)])\}\}/[u\backslash (xv)]\cdot v.
\end{displaymath}
Thus, $(f\cdot b)\theta/d=[f\cdot (a\backslash c\theta^{-1})]\theta$
if and only if identity OSI$_0^1$ is obeyed by $(Q, \cdot
,\backslash ,/)$.

The next formulae after OSI$_0^1$ derived by putting $u=v=e$ into
OSI$_0^1$. Consequently,
$T_{(x)}=\mathbb{L}_x\mathbb{R}_x\mathbb{R}_{x^\lambda}R_x$. In an
Osborn loop, $T_{(x)}=L_{x^\lambda}R_x$, so we have the DLIP.

\subsection{Application Of two Universal Osborn Loops Identities To Cryptography}
Among the few identities that have been established for universal
Osborn loops in Theorem~\ref{1:8}, we would recommend two of them;
OSI$_0^{1.1}$ and DLIP for cryptography in a similar spirit in which
the cross inverse property has been used by Keedwell \cite{phd176}.
It will be recalled that CIPLs have been found appropriate for
cryptography because of the fact that the left and right inverses
$x^\lambda$ and $x^\rho$ of an element $x$ do not coincide unlike in
left and right inverse property loops, hence this gave rise to what
is called 'cycle of inverses' or 'inverse cycles' or simply 'cycles'
i.e finite sequence of elements $x_1,x_2,\cdots ,x_n$ such that
$x_k^\rho =x_{k+1}~\bmod{n}$. The number $n$ is called the length of
the cycle. The origin of the idea of cycles can be traced back to
Artzy \cite{phd140,phd193} where he also found there existence in
WIPLs apart form CIPLs. In his two papers, he proved some results on
possibilities for the values of $n$ and for the number $m$ of cycles
of length $n$ for WIPLs and especially CIPLs. We call these "Cycle
Theorems" for now.

\paragraph{}
In Corollary~3.4 of Jaiy\'e\d ol\'a and Ad\'en\'iran \cite{phd195},
it was established that in a universal Osborn loop, $J_\lambda
=J_\rho$, 3-PAP, LSIP and RSIP are equivalent conditions.
Furthermore, in a CC-loop, the power associativity property, 3-PAPL,
$x^\rho =x^\lambda$, LSIP and RSIPL were shown to be equivalent in
Corollary~3.5. Thus, universal Osborn loops without the LSIP or RSIP
will have cycles(even long ones). This exempts groups, extra loops,
and Moufang loops but includes CC-loops, VD-loops and universal
WIPLs. Precisely speaking, non-power associative CC-loops will have
cycles. So broadly speaking, universal Osborn loops that do not have
the LSIP or RSIP or 3-PAPL or weaker forms of inverse property,
power associativity and diassociativity to mention a few, will have
cycles(even long ones). The next step now is to be able to identify
suitably chosen identities in universal Osborn loops, that will do
the job the identity $xy\cdot x^\rho =y$ or its equivalents does in
the application of CIPQ to cryptography. These identities will be
called Osborn cryptographic identities(or just cryptographic
identities).

\begin{mydef}\label{5:1}(Cryptographic Identity and Cryptographic Functional)

Let $\mathcal{Q}=(Q, \cdot ,\backslash ,/)$ be a quasigroup. An
identity $w_1(x,x_1,x_2,x_3,\cdots ,x_n )=w_2(x,x_1,x_2,x_3,\cdots ,x_n)$
where $x\in Q$ is fixed, $x_1,x_2,x_3,\cdots ,x_n\in Q$, $x\not\in
\{x_1,x_2,x_3,\cdots ,x_n\}$ is said to be a cryptographic identity(CI)
of the quasigroup $\mathcal{Q}$ if it can be written in a functional form
$xF(x_1,x_2,x_3,\cdots ,x_n )=x$ such that $F(x_1,x_2,x_3,\cdots ,x_n)\in
\mathcal{M}\textrm{ult}(\mathcal{Q})$. $F(x_1,x_2,x_3,\cdots ,x_n )=F_x$
is called the corresponding cryptographic functional(CF) of the CI
at $x$.
\end{mydef}

\begin{mylem}\label{5:2}
Let $\mathcal{Q}=(Q, \cdot ,\backslash ,/)$ be a loop with identity
element $e$ and let $CF_x(\mathcal{Q})$ be the set of all CFs in
$\mathcal{Q}$ at $x\in Q$. Then, $CF_x(\mathcal{Q})\le
\mathcal{M}\textrm{ult}(\mathcal{Q})$ and $CF_e(\mathcal{Q})\le
\textrm{Inn}(\mathcal{Q})$.
\end{mylem}
{\bf Proof}\\
The proof is easy and can be achieved by simply verifying the group axioms in $CF_x(\mathcal{Q})$ and $CF_e(\mathcal{Q})$.
\begin{enumerate}
\item
\begin{description}
\item[Closure] Obviously by definition, $CF_x(\mathcal{Q})\subset\mathcal{M}\textrm{ult}(\mathcal{Q})$. Let $F_1,F_2\in CF_x(\mathcal{Q})$. So, $xF_1F_2=xF_2=x$ which implies that $F_1F_2\in CF_x(\mathcal{Q})$.
\item[Associativity] Trivial.
\item[Identity] $xI=x$. So, $I\in CF_x(\mathcal{Q})$.
\item[Inverse] $F\in CF_x(\mathcal{Q})\Rightarrow xF=x\Rightarrow xF^{-1}=x\Rightarrow F^{-1}\in CF_x(\mathcal{Q})$.
\end{description}
$\therefore$ $CF_x(\mathcal{Q})\le
\mathcal{M}\textrm{ult}(\mathcal{Q})$.
\item Obviously by definition, $CF_e(\mathcal{Q})\subset\textrm{Inn}(\mathcal{Q})$.
The procedure of the proof that $CF_e(\mathcal{Q})\le
\textrm{Inn}(\mathcal{Q})$ is similar to that for $CF_x(\mathcal{Q})\le
\mathcal{M}\textrm{ult}(\mathcal{Q})$
\end{enumerate}

\begin{mydef}\label{5:1.1}(Degree of Cryptographic Identity and Cryptographic Functional)

Let $\mathcal{Q}=(Q, \cdot ,\backslash ,/)$ be a quasigroup and $\mathcal{I}$ an identity in $\mathcal{Q}$. If $\mathcal{I}$ is a CI with CF $F$, then the functions $F_1,F_2,F_3,\cdots F_n\in \mathcal{M}\textrm{ult}(Q)$ are called the $n$-components of $F$, written $F=(F_1,F_2,F_3,\cdots ,F_n)$ if $F=F_1\circ F_2\circ F_3\circ\cdots\circ F_n$. The maximum $n\in \mathbb{Z}^+$ such that $F=F_1\circ F_2\circ F_3\circ\cdots\circ F_n$ is called the degree of $F$ or $I$.
\end{mydef}

\begin{myexam}\label{5:1.2}
Consider a CIPQ $L$. The identity $\mathcal{I}~:~xy\cdot x^\rho =y$ is a CI at any point $y\in L$ with CF $F(x)=F_y=L_xR_{x^\rho }$. It can be seen that $F(x)=F_1(x)F_2(x)$ where $F_1(x)=L_x$ and $F_2(x)=R_{x^\rho }$, thus, $F(x)=(L_x,R_{x^\rho })$. $\mathcal{I}$ is of degree $2$. Note that an $F$ of rank $1$ is the identity mapping $I$.
\end{myexam}

\begin{mylem}\label{5:1.3}
Let $\mathcal{Q}=(Q, \cdot ,\backslash ,/)$ be a quasigroup and $\mathcal{I}$ an identity in $\mathcal{Q}$. If $\mathcal{I}$ is a CI with CF $F$ at any point $x\in Q$ such that $F=(F_1,F_2)$, then $F_1\in CF_x(\mathcal{Q})$ if and only if $F_2\in CF_x(\mathcal{Q})$.
\end{mylem}
{\bf Proof}\\
$F=(F_1,F_2)$ implies that $xF=xF_1F_2=x$. Thus, $F_1\in CF_x(\mathcal{Q})\Leftrightarrow xF_2=x\Leftrightarrow F_2\in CF_x(\mathcal{Q})$.

\begin{mylem}\label{5:1.4}
Let $\mathcal{Q}=(Q, \cdot ,\backslash ,/)$ be a quasigroup and $\mathcal{I}$ an identity in $\mathcal{Q}$. If $\mathcal{I}$ is a CI with CF $F$ at any point $x\in Q$ such that $F=(F_1,F_2,F_3,\cdots ,F_n)$, then $F_1,F_2,F_3,\cdots ,F_{n-1}\in CF_x(\mathcal{Q})$ implies $F_n\in CF_x(\mathcal{Q})$.
\end{mylem}
{\bf Proof}\\
$F=(F_1,F_2,F_3,\cdots ,F_n)$ implies that $xF=xF_1F_2F_3\cdots F_n=x$. Thus, $F_1,F_2,F_3,\cdots ,F_{n-1}\in CF_x(\mathcal{Q})\Rightarrow xF_n=x\Rightarrow F_n\in CF_x(\mathcal{Q})$.

\begin{mylem}\label{5:3}
Let $\mathcal{Q}=(Q, \cdot ,\backslash ,/)$ be a quasigroup.
\begin{enumerate}
\item $T_{(x)}\in CF_z(\mathcal{Q})$ if and only if $z\in C(x)$ for all $x,z\in Q$,
\item $R_{(x,y)}\in CF_z(\mathcal{Q})$ if and only if $z\in N_\lambda (x,y)$ for all $x,y,z\in Q$,
\item $L_{(x,y)}\in CF_z(\mathcal{Q})$ if and only if $z\in N_\rho (x,y)$ for all $x,y,z\in Q$,
\end{enumerate}
where $N_\lambda (x,y)=\{z\in Q~|~zx\cdot y=z\cdot xy\}$, $N_\rho
(x,y)=\{z\in Q~|~y\cdot xz=yx\cdot z\}$ and $C(z)=\{y\in Q~|~zy=yz\}$.
\end{mylem}
{\bf Proof}\\
\begin{enumerate}
\item $T_{(x)}\in CF_y(\mathcal{Q})\Leftrightarrow yT_{(x)}=y\Leftrightarrow yR_x=yL_x\Leftrightarrow yx=xy\Leftrightarrow y\in C(x)$.
\item $R_{(x,y)}\in CF_z(\mathcal{Q})\Leftrightarrow zR_{(x,y)}=z\Leftrightarrow zR_xR_y=zR_{xy}\Leftrightarrow zx\cdot y=z\cdot xy\Leftrightarrow z\in N_\lambda (x,y)$.
\item $L_{(x,y)}\in CF_z(\mathcal{Q})\Leftrightarrow zL_{(x,y)}=z\Leftrightarrow zL_xL_y=zL_{yx}\Leftrightarrow y\cdot xz=yx\cdot z\Leftrightarrow z\in N_\rho (x,y)$.
\end{enumerate}

\begin{mylem}\label{5:4}
Let $\mathcal{Q}=(Q, \cdot ,\backslash ,/)$ be a left universal
Osborn loop. Then, the identities OSI$_0^{1.1}$ and DLIP are CIs with degrees $6$ and $4$ respectively.
\end{mylem}
{\bf Proof}\\ From Theorem~\ref{1:8}:
\begin{description}
\item[OSI$_0^{1.1}$] is $z=x\cdot \{[x\backslash
(zx)]/x\cdot
x^\lambda\}x$, which can be put in the form $z=zR_x\mathbb{L}_x\mathbb{R}_xR_{x^\lambda}R_xL_x$. Thus, OSI$_0^{1.1}$ is a CI with CF $F(x)=R_x\mathbb{L}_x\mathbb{R}_xR_{x^\lambda}R_xL_x$ of degree $6$.
\item[DLIP] is $x^\lambda\cdot
xy)x^\lambda\cdot x=y$, which can be put in the form $yL_xL_{x^\lambda}R_{x^\lambda}R_x=y$. Thus, OSI$_0^{1.1}$ is a CI with CF $F(x)=L_xL_{x^\lambda}R_{x^\lambda}R_x$ of degree $4$.
\end{description}

\paragraph{Discussions}
Since the identities OSI$_0^{1.1}$ and DLIP have degrees $6$ and $4$  respectively, then they are "stronger" than the CIPI which has a degree of $2$ and hence
will posse more challenge for an attacker(than the CIPI) to
break into a system. As described by Keedwell, for a CIP, it is
assumed that the message to be transmitted can be represented as
single element $x$ of a CIP quasigroup and that this is enciphered
by multiplying by another element $y$ of the CIPQ so that the
encoded message is $yx$. At the receiving end, the message is
deciphered by multiplying by the inverse of $y$. But for the
identities OSI$_0^{1.1}$ and DLIP, procedures of enciphering and deciphering
are more than one in a universal Osborn loop. For instance, if the CFs of identities  OSI$_0^{1.1}$ and DLIP are $F$ and $G$, respectively such that $F=F_1F_2$ and $G=G_1G_2$ where
\begin{displaymath}
F_1=R_x\mathbb{L}_x\mathbb{R}_x,~F_2=R_{x^\lambda}R_xL_x,~G_1=L_xL_{x^\lambda}~\textrm{and}~G_2=R_{x^\lambda}R_x.
\end{displaymath}
If it is
assumed that the message to be transmitted can be represented as
single element $y$ of a universal Osborn loop and that this is enciphered
by transforming with $F_1$ or $G_1$ so that the
encoded message is $xF_1$ or $xG_1$. At the receiving end, the message is
deciphered by transforming by $F_2$ or $G_2$. Note that the components of $F$ and $G$ are not necessarily unique. This gives room for any choice of set of components. $F_1$ or $G_1$ will be called the sender's functional component(SFC) while
$F_2$ or $G_2$ will be called the receiver's functional component(RFC).

\paragraph{Many Receivers}
So far, we have considered how to secure information in a situation whereby there is just one sender and one receiver(this is the only case which the CIP is useful for). There are some other advanced and technical information dissemination patterns(which the CIP may not be applicable to) in institutions and organization such as financial institutions in which the information or data to be sent must pass through some other parties(who are not really cautious of the sensitive nature of the incoming information) before it gets to the main receiver. For instance, let us consider a network structure of an organization which has $n$ terminals. Say terminals $A_i$, $1\le i\le n$. Imagine that terminal $A_1$ wants to get a secured information across to terminal $A_n$ such that the information must pass through terminals $A_2,A_3,\cdots ,A_{n-1}$. Then, we need a CI $\mathcal{I}$ with CF $F$ of degree $n$ so that $F=(F_1,F_2,F_3,\cdots ,F_n)$. Thus, by making $F_i$ to be $A_i$'s functional component, then if the information $x$ is not to be known by $A_2,A_3,\cdots ,A_{n-1}$, we would make use of a $F$ which does not obey the hypothesis of Lemma~\ref{5:1.4}. That is, $F_1,F_2,F_3,\cdots ,F_{n-1}\not\in CF_x$. But if it is the other way round, an $F$ which obeys the hypothesis of Lemma~\ref{5:1.4} must be sort for. The advantage of a CF $F$ of higher degrees $n\ge 3$ over the CIPI relative to the number of attackers is illustrated below.
\begin{displaymath}
A_1\xrightarrow[\uparrow_{{\textrm{Attacker}~ 1}}\uparrow]{F_1~\textrm{Secured}}A_2\xrightarrow[\uparrow_{{\textrm{Attacker}~ 2}}\uparrow]{F_2~\textrm{Secured}}A_3\cdots\xrightarrow[]{}\cdots A_{n-1}\xrightarrow[\uparrow_{{\textrm{Attacker}~ n-1}}\uparrow]{F_n~\textrm{Secured}}A_n.
\end{displaymath}

\paragraph{}
Let us now illustrate with an example, the use of universal
Osborn loops for cryptography. But before then, it must be mentioned
that experts have found it very difficult to construct a
non-universal Osborn loop. According to Michael Kinyon during our
personal contact with him, there are two difficulties with using
software for looking for non-universal Osborn loops. One is that
non-Moufang, non-CC Osborn loops are very sparse: they do not start
to show up until order $16$(and the two of order $16$ happen to be
G-loops.) The other difficulty is that once you start to pass about
order $16$, the software slows down considerably. One of the two
Osborn loops that are G-loops constructed by Kinyon is shown in
Table~\ref{osbornloop1}.
\begin{myexam}
We shall now use the universal Osborn loop(it is a G-loop) of order
$16$ in Table~\ref{osbornloop1} to illustrate encoding and decoding.
\begin{description}
\item[Message:] OSBORN.
\item[CI:] DLIP.
\item[CF:] $G(x)=L_xL_{x^\lambda}R_{x^\lambda}R_x$
\item[Degree of CF:] 4.
\item[Encipherer:] $x=16,~x^\lambda =16^\lambda =10$.
\item[SFC:] $G_1=L_xL_{x^\lambda}$.
\item[RFC:] $G_2=R_{x^\lambda}R_x$.
\item[Representation($y$):] $\textrm{B}\leftrightarrow 7$, $\textrm{N}\leftrightarrow 9$, $\textrm{O}\leftrightarrow 11$, $\textrm{R}\leftrightarrow
12$, $\textrm{S}\leftrightarrow 13$.
\end{description}
The information to be transmitted is "OSBORN". The encoded
message is $(9,16,7,9,10,12)$ while the message decoded is
$(11,13,7,11,12,9)$. The computation for this is as
shown in Table~\ref{osbornloop9}.
\end{myexam}
\newpage
\begin{table}[tbp]
\begin{center}
\begin{tabular}{|c|c|c|c|}
\hline
LETTER & ENCIPHERING & DECIPHERING & DECODED LETTER \\

 & $y'=yG_1$ & $y'G_1G_2=y$ &  \\
\hline
B & $10(16\cdot 7)=7$ & $(7\cdot 10)16=7$ & $7$\\
\hline
N & $10(16\cdot 9)=12$ & $(12\cdot 10)16=9$ & $9$\\
\hline
O & $10(16\cdot 11)=9$ & $(9\cdot 10)16=11$ & $11$\\
\hline
R & $10(16\cdot 12)=10$ & $(10\cdot 10)16=12$ & $12$\\
\hline
S & $10(16\cdot 13)=16$ & $(16\cdot 10)16=13$ & $13$\\
\hline
\end{tabular}
\end{center}\caption{A Table of cryptographic Process using identity DLIP in a universal Osborn loop}\label{osbornloop9}
\end{table}
\begin{table}[!hbp]
\begin{center}
\begin{tabular}{|c||c|c|c|c|c|c|c|c|c|c|c|c|c|c|c|c|}
\hline
$\cdot $ & 1 & 2 & 3 & 4 & 5 & 6 & 7 & 8 & 9 & 10 & 11 & 12 & 13 & 14 & 15 & 16 \\
\hline  \hline
1 & 1 & 2 & 3 & 4 & 5 & 6 & 7 & 8 & 9 & 10 & 11 & 12 & 13 & 14 & 15 & 16\\
\hline
2 & 2 & 1 & 4 & 3 & 6 & 5 & 8 & 7 & 10 & 9 & 12 & 11 & 14 & 13 & 16 & 15\\
\hline
3 & 3 & 4 & 1 & 2 & 7 & 8 & 5 & 6 & 11 & 12 & 9 & 10 & 15 & 16 & 13 & 14\\
\hline
4 & 4 & 3 & 2 & 1 & 8 & 7 & 6 & 5 & 12 & 11 & 10 & 9 & 16 & 15 & 14 & 13\\
\hline
5 & 5 & 6 & 8 & 7 & 1 & 2 & 4 & 3 & 13 & 14 & 16 & 15 & 10 & 9 & 11 & 12\\
\hline
6 & 6 & 5 & 7 & 8 & 2 & 1 & 3 & 4 & 14 & 13 & 15 & 16 & 9 & 10 & 12 & 11\\
\hline
7 & 7 & 8 & 6 & 5 & 3 & 4 & 2 & 1 & 15 & 16 & 14 & 13 & 12 & 11 & 9 & 10\\
\hline
8 & 8 & 7 & 5 & 6 & 4 & 3 & 1 & 2 & 16 & 15 & 13 & 14 & 11 & 12 & 10 & 9\\
\hline
9 & 9 & 10 & 11 & 12 & 15 & 16 & 13 & 14 & 5 & 6 & 7 & 8 & 3 & 4 & 1 & 2\\
\hline
10 & 10 & 9 & 12 & 11 & 16 & 15 & 14 & 13 & 6 & 5 & 8 & 7 & 4 & 3 & 2 & 1\\
\hline
11 & 11 & 12 & 9 & 10 & 13 & 14 & 15 & 16 & 8 & 7 & 6 & 5 & 2 & 1 & 4 & 3\\
\hline
12 & 12 & 11 & 10 & 9 & 14 & 13 & 16 & 15 & 7 & 8 & 5 & 6 & 1 & 2 & 3 & 4\\
\hline
13 & 13 & 14 & 16 & 15 & 12 & 11 & 9 & 10 & 1 & 2 & 4 & 3 & 7 & 8 & 6 & 5\\
\hline
14 & 14 & 13 & 15 & 16 & 11 & 12 & 10 & 9 & 2 & 1 & 3 & 4 & 8 & 7 & 5 & 6\\
\hline
15 & 15 & 16 & 14 & 13 & 10 & 9 & 11 & 12 & 4 & 3 & 1 & 2 & 6 & 5 & 7 & 8\\
\hline
16 & 16 & 15 & 13 & 14 & 9 & 10 & 12 & 11 & 3 & 4 & 2 & 1 & 5 & 6 & 8 & 7\\
\hline
\end{tabular}
\end{center}
\caption{The first Osborn loop of order $16$ that is a
G-loop}\label{osbornloop1}
\end{table}

\end{document}